\documentclass[preprint,12pt]{elsarticle}
\usepackage[utf8]{inputenc}
\usepackage[a4paper, total={6in, 8in}]{geometry}
\usepackage{amssymb,amsmath,amsthm,latexsym}
\usepackage{color,soul}
\usepackage{xcolor}
\usepackage{tablefootnote}
\usepackage{threeparttable}
\usepackage{booktabs}
\usepackage{siunitx,caption}

\makeatletter
\renewcommand{\maketag@@@}[1]{\hbox{\m@th\normalsize\normalfont#1}}
\makeatother

\usepackage{footnote}
\usepackage{graphicx}
\usepackage{float}
\usepackage{spverbatim}
\newtheorem{thm}{Theorem}[section]

\newtheorem{lem}[thm]{Lemma}

\newtheorem{rema}[thm]{Remark}
\newtheorem{defi}[thm]{Definition}
\newtheorem{countexam}[thm]{Counterexample}
\newtheorem{examp}[thm]{Example}
{ \theoremstyle{remark} }
\journal{Journal}
\begin{document}
\begin{frontmatter}
\title{Some Results On the Primary Order Preserving  Properties of  Stochastic Orders}
\author[inst1,inst2,inst3]{Mohsen Soltanifar}
\affiliation[inst1]{organization={Analytics Division},
            addressline={College of Professional Studies}, 
            city={Northeastern University},
            postcode={Vancouver}, 
            state={BC},
            country={Canada}}
\affiliation[inst2]{organization={Real World Analytics},
            addressline={Cytel Canada Health Inc}, 
            city={Vancouver},
            state={BC},
            country={Canada}}            
\affiliation[inst3]{organization={Biostatistics Division},
            addressline={Dalla Lana School of Public Health}, 
            city={University of Toronto},
            postcode={Toronto}, 
            state={ON},
            country={Canada}}
\begin{abstract}
In this paper, we consider the problem of order preservation under addition and multiplication operators over the vector space of univariate 
real-valued random variables. Consistent with the case of usual order over the real numbers-as constant random variables, we prove that the usual stochastic order is both additive and multiplicative. We additionally discuss the situation for nine extra univariate stochastic orders presenting counterexamples in some cases. The problem remains unsolved for other univariate stochastic orders.
\end{abstract}


\end{frontmatter}

\section{Introduction}
\label{first section}
In their 1934 book on inequalities, Hardy, Littlewood and Polya introduced the concept of majorization as one of the fundamental building blocks of stochastic orders. Later in 1955, Lehmann introduced the concept of stochastic orders on real valued random variables, \cite{R-1}. Since then, inspired by their application in many fields there has been a growing literature on these orders specially from 1994. They have been playing a key role in comparison of different probability models in wide range of research areas such as survival analysis, reliability, queuing theory, biology and actuarial science. They serve as an informative comparison criteria between different distributions much more effective than low informative distributional point comparisons of means, medians, variance and IQRs. In almost all cases, they have the same fundamental properties of usual order $\leq$ in real numbers including reflexivity, anti-symmetry, transitivity.\par 
This paper deals with exploring another aspect of stochastic orders including order preserving additive and multiplicative properties. Our investigation originates from natural existence of earlier fundamental properties in both the set of real numbers equipped with usual order and the set of real-valued random variables equipped with one of stochastic orders. Given existence of  the mentioned order preserving properties on the set of real numbers equipped  with the usual order, there was a missing investigation on their parallel existential conditions for the case of the set of real-valued random variables equipped with a stochastic order.\par 
This paper is divided into four sections: preliminaries, order preserving properties of usual stochastic order; and, order preserving properties of four other stochastic orders.  In the first section, we provide the required definitions and established results for the next two sections. Then in the second section, we establish order preserving additive and multiplicative properties for the usual stochastic order. Furthermore, in the third section we discuss these properties for nine  extra univariate stochastic orders including the the moment, the Laplace transform, the increasing convex, the starshaped, the moment generating function, the convolution, the hazard rate, the likelihood ratio and the mean residual life order. We conclude the work with a discussion on the current results.\par


\section{Preliminaries}
The reader who has studied the concepts of  order and stochastic order is well acquainted with the  following definitions and results. For an essential account of the mentioned concepts, see \cite{R0,R1,R2}. We begin with some definitions:\par 
\begin{defi}
Let $P$ and $Q$ be ordered sets. A map $\phi: P\rightarrow Q$ is said to be order preserving whenever $x\leq y$ in $P$ implies $\phi(x)\leq \phi(y),$ in $Q.$ 
\end{defi}
When $P=Q$ are real valued vector spaces, and for fixed $z\in P$ the maps $\phi_{z}^{add}(x)=x+z$ and $\phi_{z}^{mul}(x)=xz (0<z)$ are order preserving, it is said that the order $\leq$ has additivity and multiplication properties. \par 
It is trivial that for the case $P=\mathbb{R}$ equipped with its usual order $\leq$, the order has additivity and multiplication properties.   From now onward, throughout this paper it is assumed that $P$ is the set of all real valued random variables. For the orders related to this $P$ we begin with one of the most well-known ones:\par 
\begin{defi}
\label{defi1} Let $X$ and $Y$ be two real-valued random variables with  associated Cumulative Distribution Function (CDF)s $F_X$, and $F_Y$, respectively. $X$ is said to be less or equal  than $Y$ in the usual stochastic order denoted by $X\leq_{st}Y$, if $F_{X}(t)\geq F_Y(t)\ \ \ \text{for all}\ \ t\in\mathbb{R}.$
\end{defi}
The following example \cite{Muller2001} has  key applications in the next section:\par 
\begin{examp}
\label{examp1}
Let $X\sim N(\mu_X,\sigma_X^2), Y\sim N(\mu_Y,\sigma_Y^2)$ be two normally distributed random variables. Then, a necessary and sufficient condition for $ X\leq_{st}Y$ is that $\mu_X\leq \mu_Y$ and $\sigma_X=\sigma_Y.$   
\end{examp}
The other orders of interest in this work are introduced in the following definition:\par   
\begin{defi}
\label{defi2}
Let $X$ and $Y$ be two real-valued random variables with  associated CDFs $F_X,F_Y$, Laplace transforms $L_X,L_Y,$ moment generating function $M_X,M_Y,$ and mean residual life functions $m_X,m_Y,$ respectively. Then, $X$ is said to be less or equal  than $Y$ in the:\newline 
(i)  moment  order, denoted by $X\leq_{m} Y$, if for $0\leq X,Y$ we have:
$ E(X^m) \leq E(Y^m)\ \  \text{for all}\ \  m\in \mathbb{N}.$\newline 
(ii)   Laplace transformation order, denoted by $X\leq_{Lt} Y$, if with assumption $0\leq X,Y$  we have: $ L_X(s) \geq L_Y(s) \ \ \text{for all}\ \  s\in\mathbb{R}^{+}.$\newline
(iii) increasing convex order denoted by $X\leq_{icx} Y,$ if :$E(g(X))\leq E(g(Y))$: for all increasing convex functions $g.$\newline 
(iv) starshaped order denoted by $X\leq_{ss} Y,$ if :$E(g(X))\leq E(g(Y))$: for all starshaped functions $g:\mathbb{R}^+\rightarrow\mathbb{R}^+.$\newline 
(v) moment generating function denoted by $X\leq_{\text{mgf}} Y$ , if with assumption $0\leq X,Y$  we have: $ M_X(s) \leq M_Y(s) \ \ \text{for all}\ \  s\in\mathbb{R}^{+}.$\newline
(vi) convolution order, denoted by $X\leq_{conv} Y$,if  for some non-negative independent random variable $U$ of $X,$ $Y=_{st} X+U,$\newline
(vii) hazard rate order, denoted by $X \leq_{hr}Y,$ if for the hazard function $r(t)=\frac{\frac{dF(t)}{dt}}{1-F(t)}$ we have:
$r_{X} (t)\geq r_{Y} (t) \ \ \text{for all}\ \  t\in\mathbb{R},$ \newline
(viii) likelihood ratio order, denoted by $X\leq_{lr}Y,$ if the likelihood ratio function $f_{ratio}(t)=\frac{\frac{dF_X(t)}{dt}}{\frac{dF_Y(t)}{dt}}$ is decreasing for all $t\in\mathbb{R}.$\newline
(ix) mean residual life order, denoted by $X\leq_{mrl}Y,$ if for the mean residual function
$m(t)=\frac{\int_{t}^{\infty}(1-F(u))du}{1-F(t)}$ we have: $m_X(t)\leq m_Y(t)$   for all $t\in\mathbb{R}.$
\end{defi}
It has been mentioned in the literature that usual stochastic order is reflexive,  anti-symmetric and transitive  \cite{R1}. Regarding order preserving maps on $P$ we have \cite{R2}:\par 
\begin{thm}
\label{thm2.4}
If $X \leq_{ord} Y$ and $\phi:\mathbb{R}\rightarrow\mathbb{R}$ is any increasing function, then $\phi(X) \leq_{ord} \phi(Y )$  where $\leq_{ord}$ refers to $\leq_{st}$, $\leq_{hr}$ or $\leq_{lr}.$
\end{thm}
In particular, for constant random variable $Z_0$,  $\phi_{Z_0}^{add}$ and $\phi_{Z_0}^{mul}$ are order preserving and hence the usual stochastic order is additive and multiplicative in this special case. However, the case in general is still unknown.\par 
An straightforward verification shows that considering the subset $P_0$ of all constant real-valued random variables of $P$ equipped with one of above four stochastic orders, the order has both of mentioned order preserving properties. This is a parallel result to the case of real numbers $\mathbb{R}$ equipped with the usual order $\leq$. The following lemma will be useful in the proof of the upcoming theorems in the next section, \cite{R1,R2}.\par 

\begin{lem}
\label{lemma1} 
Let $X,Y$ be two real valued random variables with associated CDFs $F_X,F_Y$, and Laplace transforms $L_X,L_Y,$ respectively. Then:\newline 
(i) $X\leq_{icx}Y$ if and only if:
$\int_{x}^{\infty}(1-F_X(t))dt \leq \int_{x}^{\infty}(1-F_Y(t))dt,$ for all $x\in\mathbb{R}.$\newline   
(ii) $X\leq_{ss}Y$ if and only if:
$\int_{x}^{\infty} t dF_X(t) \leq \int_{x}^{\infty}  tdF_Y(t) $ for all $x\in\mathbb{R}_{0}^{+}.$\newline (iii) $X\leq_{conv}Y$ if and only if  for $\phi_{X,Y}=\frac{L_Y}{L_X}$ we have: $(-1)^n.\phi_{X,Y}^{(n)}(s) \geq 0\ \ \ \text{for all}\ \  (s\in\mathbb{R}^{+}, \ n\in \mathbb{N}). $\newline
(iv) $X\leq_{mgf}Y$ if and only if: $-Y\leq_{Lt}-X.$\newline
(v) $X\leq_{mrl}Y$ with given finite means if and only if $g_{ratio}(t)=\frac{\int_{t}^{\infty}(1-F_X(u))du}{\int_{t}^{\infty}(1-F_Y(u))du}$ is  decreasing.   
\end{lem}

We conclude this section with a remark on the independency of involved random variables:

\begin{rema}
\label{rem1}
Considering any pair of real numbers $x_0,y_0$ as constant random variables $X=x_0$ and $Y=y_0$, it is trivial that any constant random variable $Z=z_0$ is independent from each of them. Hence, it is natural to maintain this independency assumption for inferential results for the case of non-constant random variables. 
\end{rema}
\section{Order Preserving Properties of Usual Stochastic Order}
\label{math}
This section deals with general properties of usual stochastic order. Throughout this section and the next one, we assume the given random variable  $Z$ is independent from $X$ and $Y$. The first three properties has been established as mentioned above. We discuss the later two properties and present an evidence of the necessity of the independence condition mentioned in Remark \ref{rem1}.\par 

\begin{thm}\label{thm3.1}
The usual stochastic order is (i) additive; and, (ii) multiplicative.
\end{thm}
\textbf{Proof.} This is a direct consequence from \cite{R2} Theorem 1.A.3(b) with $m=2, X_1=X, X_2=Z, Y_1=Y, Y_2=Z,$ and consideration of two bi-variate increasing functions $\psi_1(u,v)=u+v$ (for additivity) and $\psi_2(u,v)=uv$ (for multiplication). $\Box$\newline\\ 
An immediate application of the Theorem \ref{thm3.1} is in the area of computational psychology for comparing reaction times between participants. Here, for the case of the well-known Ex-Gaussian(ExG) distribution \cite{Heatcote1996} defined by $ExG(\mu,\sigma,\tau)\overset{d}{=}N(\mu,\sigma^2)+Exp(\tau)$ we may characterize which experimental participant is faster than the other in terms of the fitted distributional parameters:\par  
\begin{examp}
\label{examp2}
Let $ExG(\mu_X,\sigma_X,\tau_X),ExG(\mu_Y,\sigma_Y,\tau_Y)$ be two reaction times Ex-Gaussian distributions.  Then, by an application of Example \ref{examp1}, two applications of Theorem \ref{thm3.1}, and considering the transitivity of the usual stochastic order, a sufficient condition for $ExG(\mu_X,\sigma_X,\tau_X) \leq_{st} ExG(\mu_Y,\sigma_Y,\tau_Y)$ is $\mu_X\leq \mu_Y$,  $\sigma_X=\sigma_Y,$ and   $\tau_X\geq \tau_Y.$ 
\end{examp}
  
\begin{countexam}
\label{countexam1}
The assumption of independence is necessary for the usual stochastic order additivity. As a counterexample, let $X\sim N(0,\frac{1}{2}), Y\sim N(1,\frac{1}{2})$ be two independent normal random variables. Then, by Example \ref{examp1} for $\mu_X=0<1=\mu_Y, \sigma_X=\frac{1}{2}=\sigma_Y$ we have: $X\leq_{st}Y.$ On the other hand, take $Z=-X$, consider $Y-X \sim N(1,1),$ and $F_{0}(-1)=1_{[0,\infty)}(-1)=0 < \Phi(-2)=F_{Y-X}(-1).$ Consequently: $0\nleq_{st}Y-X.$ Figure \ref{FigA1} (a) presents this case.
\end{countexam}

\begin{countexam}
\label{countexam11}
The assumption of independence is necessary for the usual stochastic order multiplication. As a  counterexample, let $X\sim logN(0,\frac{1}{2}) , Y\sim logN(1,\frac{1}{2})$ be two independent lognormal random variables with $X\leq_{st}Y,$ (this is a consequence from Example \ref{examp1} for $\mu_X=0<1=\mu_Y,\sigma_X=\frac{1}{2}=\sigma_Y$ and Theorem \ref{thm2.4} for $\phi_1(u)=\exp(u)$). Consider dependent random variable $Z=X^{-1}$ (this is a consequence of strict Jensen's inequality). If $XZ\leq_{st}YZ,$ then by Theorem \ref{thm2.4} (for $\phi_2(u)=\log(u)$) we have $0\leq_{st} N(1,1),$ a contradiction. Hence, $XZ\nleq_{st}YZ.$ 
\end{countexam}

\begin{countexam}
\label{countexam111}
The assumption of positivity of random variable $Z$ is necessary for the usual stochastic order multiplication. As a counterexample, let $X_{\mu,\sigma}\sim N(\mu,\sigma^2)$ and $Z\sim 0.5(\delta_{\{-1\}}+\delta_{\{+1\}}).$ Then, $F_{X_{\mu,\sigma}Z}(t)=0.5(\Phi((t-\mu)/\sigma)+\Phi((t+\mu)/\sigma))$ for all $t\in\mathbb{R}.$ Take $X=X_{0,1}$ and $Y=X_{1,1},$ then,  $X\leq_{st}Y$ (this is a consequence from Example \ref{examp1} for $\mu_X=0<1=\mu_Y,\sigma_X=1=\sigma_Y$). But, given $F_{XZ}(-1)=0.1586553<0.2613751=F_{YZ}(-1)$, we have $XZ\nleq_{st}YZ.$ Figure \ref{FigA1} (b) presents this case.
\end{countexam}
\section{Order Preserving Properties of Other Stochastic Orders}
\label{math}
As we mentioned, there are other types of stochastic orders in the literature of particular interest. In this section we deal with the similar problem of determining their main  properties as the case of usual stochastic order for the moment, the Laplace transform, the increasing convex, the starshaped, the moment generating function, the convolution, the hazard rate, the likelihood ratio and the mean residual life order. First of all, for the case of reflexivity, anti-symmetry and transitivity, it is easy to show that all of these orders are preserving these properties. Next, we deal with the additivity and multiplication as follows:\par  
 
\begin{thm}
\label{thm4.1}
Consider the moment order, the Laplace transformation order, the increasing convex order, the starshaped order, the moment generating function order, and the convolution. Then: (i) All six orders are additive; (ii) The first five orders are multiplicative as well.
\end{thm}
\textbf{Proof.}
To prove(i), First, for the moment order, the proof is straightforward application of binomial identity, and the independence condition. Indeed,  let $X\leq_{m}Y$ and $X,Y>0$. Then, for any independent $0\leq_{m}Z$ we have:
\begin{eqnarray}
E((X+Z)^m)&=& E(\sum_{k=0}^{m} C_k^m X^k Z^{m-k})= \sum_{k=0}^{m} C_k^m E(X^k)E(Z^{m-k}) 
           \leq \sum_{k=0}^{m} C_k^m E(Y^k)E(Z^{m-k})\nonumber\\
           &=&E(\sum_{k=0}^{m} C_k^m Y^k Z^{m-k})
          = E((Y+Z)^m)\ \ \forall m\in\mathbb{N}.\nonumber
\end{eqnarray} 
Thus, $X+Z\leq_{m} Y+Z.$\newline 
Second, for the Laplace transformation, the proof is straightforward from of the Laplace transformation multiplication of sum of independent random variables. Indeed, let $X\leq_{Lt}Y$ and  $0<_{Lt}Z$ be independent. Then we have:
\begin{eqnarray}
L_{X+Z}(t)&=&E(e^{-t(X+Z)})=E(e^{-tX}e^{-tZ})=E(e^{-tX})E(e^{-tZ})\leq E(e^{-tY})E(e^{-tZ})\nonumber\\
&=&E(e^{-tY}e^{-tZ})=E(e^{-t(Y+Z)}) =L_{Y+Z}(t), \forall t\in \mathbb{R},\nonumber
\end{eqnarray}
implying  $X+Z\leq_{Lt}Y+Z.$ \newline 
Third, for the increasing convex order, let $X\leq_{icx}Y$ and $Z$ be independent from $X,Y.$ Then by Lemma \ref{lemma1} and two times application of Fubini's theorem it follows that:
\begin{eqnarray}
\int_{x}^{\infty}(1-F_{X+Z}(t))dt&=&\int_{x}^{\infty}(1- \int_{z=-\infty}^{\infty}F_X(t-z)dF_Z(z))dt=\int_{x}^{\infty} \int_{z=-\infty}^{\infty} (1-F_X(t-z))dF_Z(z)dt\nonumber\\ 
                                 &=& \int_{z=-\infty}^{\infty} \int_{x}^{\infty}(1-F_X(t-z))dt dF_Z(z) 
                                 \leq  \int_{z=-\infty}^{\infty} \int_{x}^{\infty}(1-F_Y(t-z))dtdF_Z(z)\nonumber\\ &=& \int_{x}^{\infty} \int_{z=-\infty}^{\infty} (1-F_Y(t-z))dF_Z(z)dt=\int_{x}^{\infty}(1- \int_{z=-\infty}^{\infty}F_Y(t-z)dF_Z(z))dt\nonumber\\
                                 &=& \int_{x}^{\infty}(1-F_{Y+Z}(t))dt\ \ \text{for all}\ \ x\in\mathbb{R}.\nonumber 
\end{eqnarray} 
Now, by another application of Lemma \ref{lemma1} it follows that $X+Z\leq_{icx}Y+Z.$\newline\\
Fourth, the starshaped order, let $X\leq_{ss}Y$ and $Z$ be independent from $X,Y.$ Then by Lemma \ref{lemma1} and two times application of Fubini's theorem it follows that:
\begin{eqnarray}
\int_{x}^{\infty}tdF_{X+Z}(t)&=&
\int_{x}^{\infty}td (\int_{z=-\infty}^{\infty}F_Z(z)dF_X(t-z)) 
=
\int_{x}^{\infty}  \int_{z=-\infty}^{\infty}tdF_Z(z)dF_X(t-z) 
\nonumber\\
&=&
\int_{z=-\infty}^{\infty} (\int_{x}^{\infty}  tdF_X(t-z))dF_Z(z) 
\leq
\int_{z=-\infty}^{\infty} (\int_{x}^{\infty}  tdF_Y(t-z))dF_Z(z) 
\nonumber\\
&=&
\int_{x}^{\infty}  \int_{z=-\infty}^{\infty}tdF_Z(z)dF_Y(t-z) 
=
\int_{x}^{\infty}td (\int_{z=-\infty}^{\infty}F_Z(z)dF_Y(t-z)) 
\nonumber\\
&=&\int_{x}^{\infty}tdF_{Y+Z}(t)
 \ \ \text{for all}\ \ x\in\mathbb{R}.\nonumber 
\end{eqnarray} 
Again, by another application of Lemma \ref{lemma1} it follows that $X+Z\leq_{ss}Y+Z.$\newline\\
Fifth, this is direct consequence from the additivity for Laplace transformation order and Lemma \ref{lemma1}.\newline 
Finally, we consider the convolution order. Let $X\leq_{conv}Y$ and $Z$ be independent from $X,Y.$ Then, by the later condition it follows:
$$\phi_{X+Z,Y+Z}(s)=\frac{L_{Y+Z}(s)}{L_{X+Z}(s)}=\frac{L_{Y}(s)L_{Z}(s)}{L_{X}(s)L_{Z}(s)}=\phi_{X,Y}(s),\ \  \text{for all}\ s\in\mathbb{R}^{+} .$$
Now, two applications of Lemma \ref{lemma1} and former condition implies $X+Z\leq_{conv}Y+Z.$ \newline 
To prove (ii), first of all, for the moment order the proof is trivial from independence condition. Second, for the Laplace transformation order, let $X\leq_{Lt}Y$ and $0<_{Lt}Z$ be independent from $X,Y.$ Then, we have:
\begin{eqnarray}
L_{XZ}(s)=\int_{z=0}^{\infty}L_{X}(sz)dF_{Z}(z)\geq \int_{z=0}^{\infty}L_{Y}(sz)dF_{Z}(z)&=&L_{YZ}(s)\ \  \text{for all}\  s\in\mathbb{R}^{+}.\nonumber
\end{eqnarray}
Hence, $XZ\leq_{Lt}YZ$.\newline
Third, and fourth, for the increasing convex order and the starshaped order,  the proof is similar to the presented proof for the additivity case above with $"X+Z", "Y+Z", "t-z"$ replaced by $"XZ", "YZ", "t/z",$ respectively.\newline 
Finally, for the moment generating function order, the proof is similar to the proof of additivity presented above. $\Box$ 
\begin{rema}
\label{rem2}
In Theorem \ref{thm4.1} the case of additivity and multipliciation for the increasing convex order can be proved similarly to the case usual order presented in the proof of \ref{thm3.1} using Theorem 4.A.15-  \cite{R2} with $g_1(u,v)=u+v;  g_2(u,v)=u*v,$ respectively. 
\end{rema}

Establishing preservation results for the orders mentioned in Theorem \ref{thm4.1}, we now discuss the cases in which the additivity and multiplicity are not preserved. We begin with the case for convolution order:\par
\begin{figure}[H] 
\begin{center}
\includegraphics[totalheight=15.5 cm, width=14 cm]{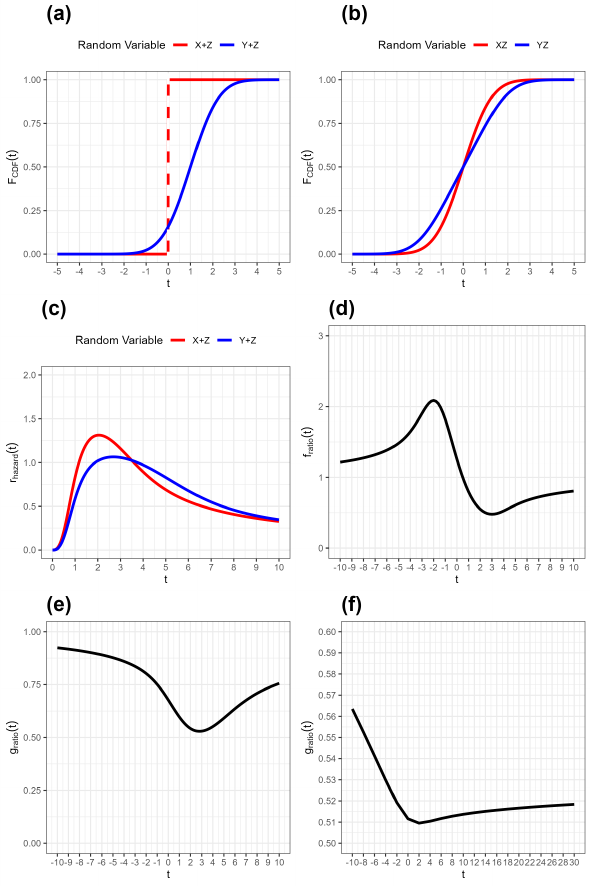}
\caption{Plots for Counterexamples: (a) Counterexample 3.3; (b) Counterexample 3.5; (c) Counterexample 4.4; (d) Counterexample 4.6; (e) Counterexample 4.7; (f) Counterexample 4.8.\label{FigA1}}\label{solfig}
\end{center}
\end{figure}

  
\begin{countexam}
\label{countexam5} 
The convolution order is not multiplicative. As the counterexample, let $X \sim Exp(1), Y\sim Exp(0.5),$ and $Z\sim Bernoulli(0.5),$ be independent from $X,Y.$ Since $\phi_{X,Y}(s)=\frac{s+1}{2s+1},$ for all $s\in\mathbb{R}^{+}$ and $(-1)^{n}\phi_{X,Y}^{(n)}(s)=\frac{n!.2^{n-1}}{(2s+1)^{n+1}}>0,$ for all $s\in\mathbb{R}^{+};$ by an application of Lemma \ref{lemma1}, $X\leq_{conv}Y.$ However, a simple conditioning on $Z$ yields   $\phi_{XZ,YZ}(s)=\frac{s^2+2s+1}{s^2+2.5s+1},$ for all $s\in\mathbb{R}^{+}.$ But, $\phi_{XZ,YZ}^{'}(s)=\frac{0.5(s^2-1)}{(s^2+2.5s+1)^2},$ changes signs from $0<s<1$ to $1\leq s.$ Consequently, by another application of Lemma \ref{lemma1}, it follows that $XZ\nleq_{conv}YZ.$
\end{countexam}
Although, for the hazard rate order $\leq_{hr}$ an special case of additive property has been established (whenever $Z$ is IFR  i.e., $r_{Z}$ is increasing),  \cite{R2} (Lemma 1.B.3), the general case does not hold as the following counterexample explains:
\begin{countexam}
\label{countexam3}
The hazard rate order is not additive. As the counterexample, let $X \sim Exp(2.5), Y\sim Exp(1.5),$ and $Z\sim loglogistic(shape=3,scale=1),$ be independent from $X,Y.$ Then, as $r_X(t)=2.5\geq 1.5=r_Y(t)\ \ \ all\ \ t\in\mathbb{R},$ it follows that $X\leq_{hr}Y.$ However, given that $r_{X+Z}(5)=0.6887492<0.8289720=r_{Y+Z}(5)$ (See Appendix B), it follows that: $X+Z\nleq_{hr}Y+Z.$ Figure \ref{FigA1} (c) presents this case.
\end{countexam}
Next, for the case of multiplicity of the hazard rate order $\leq_{hr},$ the special case for the IFR random variable $Z$ can be proved similarly to the case of additivity. However, this order does not preserve the multiplicity as the following counterexample shows:
\begin{countexam}
\label{countexam4} The hazard rate order is not multiplicative. As a counterexample, let $X^*= \exp(X), Y^*=\exp(Y),$ and $Z^*=\exp(Z)$ where $X,Y,$ and $Z$ are defined as in counterexample \ref{countexam3}. Then, $X^*\leq_{hr}Y^*$ (using Theorem \ref{thm2.4} for $\phi(u)=\exp(u)$). However, if $X^*Z^*\leq_{hr}Y^*Z^*,$ then $X+Z\leq_{hr}Y+Z$ (using Theorem \ref{thm2.4} for $\phi(u)=\log(u)$), yielding a contradiction. Hence, $X^*Z^*\nleq_{hr}Y^*Z^*.$
\end{countexam}
Furthermore, for the likelihood ratio order $\leq_{lr}$ an special case of additive property has been established (whenever $X,Y,Z$ all have logconcave densities),  \cite{R2} (Theorem.1.C.9), the general case does not hold as the following counterexample explains:
\begin{countexam}
\label{countexam5}
The likelihood ratio order is not additive. As the counterexample, let $X \sim N(0,1), Y\sim N(1,1),$ and $Z\sim Cauchy(location=0,scale=1),$ be independent from $X,Y.$ Then, it is trivial that $X\leq_{lr}Y.$ However, given that given that for  $t_1=1<t_2=3<t_3=7,$ we have $f_{ratio}(t_1)=0.7996384>0.4814513=f_{ratio}(t_2)<0.7214805= f_{ratio}(t_3),$ (See Appendix C) it follows that  it follows that: $X+Z\nleq_{lr}Y+Z.$ Figure \ref{FigA1} (d) presents this case.
\end{countexam}
Next, using the same idea of exponentiating and back and forth usage of Theorem \ref{thm2.4} described in Counterexample \ref{countexam4}, it is easily proved that the likelihood ratio order $\leq_{lr}$ is not multiplicative either.\par 
We conclude this section with the case of mean residual life order. Here, similar to the hazard rate order  an special case of additive property has been established (whenever $Z$ is IFR),  \cite{R1} (Theorem.2.4.22), the general case does not hold as the following counterexample explains it:\par 
\begin{countexam}
\label{countexam6}
The mean residual life order is not additive. As the counterexample, let $X \sim N(0,1), Y\sim N(1,2),$ and $Z\sim logNormal(meanlog=0.25,sdlog=1),$ be independent from $X,Y.$ Then, by \cite{R1}(Example.2.4.16, page 55) it follows that $X\leq_{mrl}Y.$ However, given that given that for $t_1=1<t_2=3<t_3=5,$ we have $g_{ratio}(t_1)=0.5973625>0.5301898=g_{ratio}(t_2)<0.5913812= g_{ratio}(t_3),$ (See Appendix D) it follows that  it follows that: $X+Z\nleq_{mrl}Y+Z.$ Figure \ref{FigA1} (e) presents this case.
\end{countexam}
Finally,  and a variation of counterexample \ref{countexam6} presented in the following counterexample establish the case for multiplication:\par 
\begin{countexam}
\label{countexam7}
The mean residual life order is not multiplicative. As the counterexample, let $X \sim N(0,1), Y\sim N(1,2),$ and $Z=exp(Z0): Z_0\sim logNormal(meanlog=0.25,sdlog=1) ,$ be independent from $X,Y.$ Then, $X\leq_{mrl}Y.$ However, given that given that for $t_1=-10<t_2=2<t_3=30,$ we have $g_{ratio}(t_1)=0.5632457>0.5092787=g_{ratio}(t_2)<.0.5181358= g_{ratio}(t_3),$ (See Appendix E) it follows that  it follows that: $XZ\nleq_{mrl}YZ.$ Figure \ref{FigA1} (f) presents this case.
\end{countexam}

\section{Discussion}
\subsection{Summary}
This work presented an introductory report on fundamental order preserving properties of addition and multiplication for the case of real-valued univariate random variables. It generalized the case from the usual order on the real numbers to the usual stochastic order maintaining the assumptions of independence of involved random variables $(X,Z)$ and $(Y,Z)$ and the positivity of the random variable $Z.$ However, investigation over other types of univariate stochastic orders showed that preservation of these fundamental properties varied by the type of considered order. Table \ref{table1} summarizes the main results in both sections 3 and 4:\par

\begin{table}[H]
\caption{A Summary of order preserving properties of ten univariate stochastic orders\label{table1}} 
\begin{threeparttable}[b]
\centering 
\begin{tabular}{ll ccccc} 
\hline 
& &\multicolumn{5}{c}{Order Preserving Property} \\ [0.5ex]
\hline 
\# & Order Type &  &  & & Addition& Multiplication\\[1ex] 
\hline 
i&  Usual $\leq_{st}$   & &  &  & $\checkmark$ & $\checkmark$ \\[1ex] 
ii &Moment $\leq_{m}$     &  &  & &  $\checkmark$ & $\checkmark$  \\[1ex] 
iii &Laplace transformation $\leq_{Lt}$ & & & &$\checkmark$ & $\checkmark$ \\[1ex]
iv & Increasing Convex $\leq_{icx}$   &  & & &  $\checkmark$ & $\checkmark$ \\[1ex] 
v & Starshaped $\leq_{ss}$ & & & & $\checkmark$ & $\checkmark$ \\[1ex]
vi  & Moment generating function $\leq_{mgf}$& & & &  $\checkmark$&  $\checkmark$\\[1ex]
vii & Convolution $\leq_{conv}$ &  &  &  & $\checkmark$  & $X$ \\[1ex]
viii&Hazard Rate $\leq_{hr}$  &  & &  & $X$   &$X$ \\[1ex] 
ix & Likelihood ratio $\leq_{lr}$  & & & & $X$  & $X$  \\[1ex]
x & Mean residual life $\leq_{mrl}$  & & & & $X$  & $X$  \\[1ex]
\hline 
\end{tabular}
\end{threeparttable}
\end{table}

\subsection{Limitations}
The limitations of this work are clear. First, there are other types of univariate stochastic orders that this study did not cover their preservation properties. Some of them included the  harmonic mean residual life order, Lorenz order, dilation order, dispersive order, execessive wealth order, peakedness order, pth order, star order, superadditive order, factorial moment order, and total time on test order,\cite{R1,R2}. Second,  the case for the multivatiate random variables is still unknown. Third, we saw limited applications of the preservation of the orders (as in the case of usual stochastic order) and it is plausible to investigate them furthermore. Finally, while we investigated the preservation of the properties for given stochastic order, it is plausible to identify a necessary (and sufficient) condition for a given stochastic order for preservation of this fundamental properties.\par 

\subsection{Conclusion}
This work has presented  partial solution to the problem of order preservation under summation and multiplication for the case of univariate real-valued random variables. It showed that depending to the order under consideration preservation status of these fundamental properties varies clearing the way for more investigation in the future.\par 
$$$$
\appendix
\noindent \textbf{Appendix}\newline

\section[\appendixname~\thesection]{Plots for Counterexamples}

R software code \cite{R2021} for computing the counterexamples presented in the work:\newline\\ 
\noindent \textbf{Appendix B: hazard rate order (additivity)}
\begin{spverbatim} 
#Step(1): Call Libraries in R
myRpackages<-c("distr","flexsurv","distrEx")
lapply(myRpackages, require, character.only = TRUE)
#Step(2): Compute the hazard of sum  of loglogistic(shape = a, scale = b)
and  exponential(rate=c) distributions
h_sum=function(a,b,c,t){ U<- AbscontDistribution(d = function (t) dllogis(t, shape = a, scale = b, log = FALSE), withStand = TRUE); V<- Exp(rate=c); W<- U+V;
hazard_t= (d(W)(t))/(1-p(W)(t)); return(hazard_t)}
#Step(3): Compute Counterexamples:  
h_sum(3,1,2.5,5)  #>[1] 0.6887492
h_sum(3,1,1.5,5)  #>[1] 0.8289720
\end{spverbatim}
$$$$
\noindent \textbf{Appendix C: likelihood ratio order (additivity) }
\begin{spverbatim} 
#Step(1): Call Libraries in R
myRpackages<-c("distr","flexsurv","distrEx")
lapply(myRpackages, require, character.only = TRUE)
#Step(2): Call the random variables X=N(0,1), Y=N(1,1), Z=Cauchy(location=0,
scale=1) distroptions("DefaultNrFFTGridPointsExponent" = 10)
X<- Norm(0,1)
Y<- Norm(1,1)
Z<- AbscontDistribution(d = function (t)  1/(pi*(1+t^2)), withStand = TRUE)
#Step(3): Compute the density of random variables X+Z, and Y+Z and the ratio function
XsumZ <-  convpow(X+Z,1)
YsumZ <-  convpow(Y+Z,1)
f_XsumZ=function(t){ d(XsumZ)(t)}
f_YsumZ=function(t){ d(YsumZ)(t)}
f_ratio=function(t){f_XsumZ(t)/f_YsumZ(t)}
#Step(4): Compute Counterexamples:  
f_ratio(1)  #[1] 0.7996384
f_ratio(3)  #[1] 0.4814513
f_ratio(7)  #[1] 0.7214805
\end{spverbatim}
$$$$
\noindent \textbf{Appendix D: mean residual life order (additivity) }
\begin{spverbatim} 
#Step(1): Call Libraries in R
myRpackages<-c("distr","flexsurv","distrEx","survival")
lapply(myRpackages, require, character.only = TRUE)
#Step(2): Call the random variables X=N(0,1), Y=N(1,2),
Z=logNormal(meanlog=0.25, sdlog=1)
distroptions("DefaultNrFFTGridPointsExponent" = 10)
X<- Norm(0,1)
Y<- Norm(1,2)
Z<-Lnorm(meanlog=0.25,sdlog=1) 
#Step(3): Compute the  the ratio function g
g_ratio<-function(t,u){
  S_XsumZ<- function(x) {1-p((X+Z))(x)};
  Answer1<-integrate(S_XsumZ, lower = t, upper = u)$value;
  S_YsumZ<- function(x) {1-p((Y+Z))(x)};
  Answer2<-integrate(S_YsumZ, lower = t, upper = u)$value;
  Answer=as.numeric(Answer1)/as.numeric(Answer2);
  return(as.numeric(Answer))
}
#Step(4): Compute Counterexamples:  
g_ratio(1,Inf)  #[1] 0.5973625
g_ratio(3,Inf)  #[1] 0.5301898
g_ratio(5,Inf)  #[1] 0.5913812
\end{spverbatim}
$$$$
\noindent \textbf{Appendix E: mean residual life order (multiplication) }
\begin{spverbatim} 
#Step(1): Call Libraries in R
myRpackages<-c("distr","flexsurv","distrEx","survival")
lapply(myRpackages, require, character.only = TRUE)
#Step(2): Call the random variables X=N(0,1), Y=N(1,2),
Z=logNormal(meanlog=0.25, sdlog=1)
distroptions("DefaultNrFFTGridPointsExponent" = 10)
X<- Norm(0,1)
Y<- Norm(1,2)
Z0<-Lnorm(meanlog=0.25,sdlog=1) 
#Step(3): Compute the  the ratio function g
g_ratio<-function(t,u){
  S_XmulZ<- function(x) {1-p((X*exp(Z0)))(x)};
  Answer3<-integrate(S_XmulZ, lower = t, upper = u)$value;
  S_YmulZ<- function(x) {1-p((Y*exp(Z0)))(x)};
  Answer4<-integrate(S_YmulZ, lower = t, upper = u)$value;
  Answer=as.numeric(Answer3)/as.numeric(Answer4);
  return(as.numeric(Answer))
}
#Step(4): Compute Counterexamples:  
g_ratio(-10,1000)  #[1] 0.5632457
g_ratio(2,1000)    #[1] 0.5092787
g_ratio(30,1000)   #[1] 0.5181358
\end{spverbatim}


\end{document}